\documentclass{article}
\usepackage{hyperref}
\usepackage{amsmath,amssymb, amsmath}

\usepackage{cite}
\usepackage{graphicx}
  \usepackage{hyperref}
  \usepackage{tikz}
  \usepackage{flowchart}
\usetikzlibrary{calc}
\usetikzlibrary{decorations.pathmorphing,patterns}
\usetikzlibrary{calc,patterns,decorations.markings}
\usetikzlibrary{decorations.pathreplacing}
\usetikzlibrary{shapes,arrows,spy,positioning,snakes}
\usepackage{adjustbox}
\usepackage[margin=1cm]{caption}

\usepackage{fancyhdr}

\usepackage[utf8]{inputenc}
\usepackage[T1]{fontenc}
\usepackage{etex}
\usepackage{caption}
\newcommand{\RNum}[1]{\uppercase\expandafter{\romannumeral #1\relax}}
\newcommand{\meth}[1]{\texttt{\textbf{#1}}}

\newcommand{\bigO}{\mathcal{O}}
\usepackage{geometry}

\providecommand{\keywords}[1]{\textit{Keywords:} #1}

\begin{document}

\title{An Optimization Approach to the Ordering Phase
of an Attended Home Delivery Service}

\author{
G. Cwioro
\thanks{Department of   Mathematics, Alpen-Adria Universität Klagenfurt, Austria,
 \href{mailto:guenther.cwioro@aau.at}{guenther.cwioro@aau.at}
 }
  ,
P.\ Hungerländer
\thanks{
Laboratory for Information \& Decision Systems, Massachusetts Institute of Technology
Cambridge, MA 02139, USA,
  \href{mailto:philipp.hungerlaender@aau.at}{philipp.hungerlaender@aau.at}
  }
   ,
  K.\  Maier
  \thanks{Department of Mathematics, Alpen-Adria Universität Klagenfurt, Austria,
   \href{mailto:kerstin.maier@aau.at}{kerstin.maier@aau.at}
  }
    ,
   J.\ Pöcher
   \thanks{Department of   Mathematics, Alpen-Adria Universität Klagenfurt, Austria,
    \href{mailto:joerg.poecher@aau.at}{joerg.poecher@aau.at}
    }
    , and
  C.\  Truden
  \thanks{Department of Mathematics, Alpen-Adria Universität Klagenfurt, Austria,
   \href{mailto:christian.truden@aau.at}{christian.truden@aau.at}
   }
   }

\maketitle

\abstract{
Attended Home Delivery (AHD) systems are used
whenever a supplying company offers online shopping services
that require that customers must be present when their
deliveries arrive. Therefore, the supplying company and the
customer must both agree on a time window, which ideally is
rather short, during which delivery is guaranteed. Typically, a
capacitated Vehicle Routing Problem with Time Windows forms
the underlying optimization problem of the AHD system. In this
work, we consider an AHD system that runs the online grocery
shopping service of an international grocery retailer.

The ordering phase, during which customers place their
orders through the web service, is the computationally most
challenging part of the AHD system. The delivery schedule
must be built dynamically as new orders are placed. We
propose a solution approach that
allows to
(non-stochastically) determine which delivery time windows can
be offered to potential customers. We split the computations
of the ordering phase into four key steps. For performing
these basic steps we suggest both a heuristic approach and
a hybrid approach employing mixed-integer linear programs.
In an experimental evaluation we demonstrate the efficiency of
our approaches.
}


\maketitle

\section{Introduction}
\label{sec:intro}

In recent years, online grocery shopping has gained increased popularity in
several  countries,
such as the United Kingdom where about 6.3\% \cite{Syndy2015} of all grocery shopping is bought online.
Nowadays, all major supermarket chains provide online shopping services, where customers
select groceries as well as a delivery time window on the supermarket's website.
This provides several benefits to the customers, such as 24-hour opening hours
of the online store,
quicker shopping times, the avoidance of traveling times, no carrying of heavy
or bulky items and
facilitated access for citizens with reduced mobility.
Despite of the benefits for the customers, e-grocery shopping services pose several
interrelated logistic and optimization challenges to the supplying companies.
Especially the \textit{Ordering Phase}, during which customers place their
orders, imposes a computationally challenging problem.

In this paper, we tackle this challenge in the context of a large international
supermarket chain that offers online grocery shopping.
E-grocery services are a paradigm for Attended Home Delivery Problems (AHD)
\cite{camsa05,agca08, Ehmke2012, Ehmke2014}
where the customers must be present for their deliveries. In order to ensure customer
satisfaction and to minimize undeliverable orders, it is crucial that the
supplying company provides a wide selection of rather narrow delivery time windows.
Hence, in this work, we aim to provide a framework that (non-stochastically)
determines the available time windows and dynamically builds the delivery schedule
during the \textit{Ordering Phase}.

This paper is organized as follows. First
we provide an
overview of the logistic process behind the considered AHD system and
 discuss the \textit{Ordering Phase}
in detail.
In Section \ref{sec:algorithms} we introduce the related optimization problem
and suggest algorithmic strategies for solving it.
In Section \ref{sec:experiments} we demonstrate the efficiency of our
solution approaches on benchmark instances related to an online grocery shopping
system. Finally, Section \ref{sec:conclusion} concludes the paper.

\subsection{The Attended Home Delivery Process}
\label{subsec:adhp}

Let us start with giving a short overview of the overall planning and
fulfillment
process behind the \textit{Attended Home Delivery} service,
where we describe the actions taken by the supplying e-grocery
retailer in order to fulfill the deliveries of a single day.

\paragraph{Tactical Planning Phase - (several months or weeks before delivery):}

\begin{itemize}
  \item A fleet of vehicles is set up and operation times of those vehicles are defined.
  \item Drivers are assigned to the vehicles in accordance to the legal
  regulations concerning drive and rest times.
  \item
  The supplying company defines the set of possible delivery time windows
   that will be offered to the customers through the web service.
  \end{itemize}

\paragraph{Ordering Phase - (several weeks up to days/hours before delivery):}
This phase begins once the web service starts to allow booking of
delivery time windows for the specific day of delivery.
Hence, the system must handle the following tasks:
\begin{itemize}
  \item Customers use the web service through a web site or a mobile app
  to place their orders.
  \item The system must decide which delivery time windows can be offered to a
  specific customer
  such that the delivery can be fulfilled within the time window.
  Only the resources that have been assigned during the tactical planning phase
  are available during the ordering phase.
  \item
  Once a customer has booked a delivery time window, the system must adapt the
  existing delivery schedule to accommodate the respective order.
  Furthermore, the system periodically tries to improve the
  current schedule.
\end{itemize}
During this phase, the objective is to accept as many customers as possible,
while offering as many time windows as possible to each potential customer in order to
 achieve a high degree of customer satisfaction and to also ensure good resource
utilization which translates to the overall logistics operations being cost efficient.

\paragraph{Preparation Phase - (days/hours before delivery):}

This phase is triggered once the system does no longer accept new orders through
the web service.
The objective function is now changed to minimization of the transportation costs
(overall fulfillment costs).
Another relevant aspect to be considered is the traffic flow
at the depot and the vehicle loading bays.
Hence, the system must handle the following tasks:
\begin{itemize}
  \item The delivery schedule is improved regarding the
  new objective function.
  \item Meanwhile, at the depot, the ordered goods are fetched from
  storage and consolidated accordingly to the customer orders.
\end{itemize}

\paragraph{Delivery Phase:}
In this phase the vehicles are first packed with the consolidated orders and
  prepared to leave the depot. Then the vehicles visit the customers
  according to the delivery schedule, which was generated by the
  system, such that the customers receive their orders within the selected time windows.

\subsection{Related work}

First, let us  give a short overview of related work by the authors:

\begin{itemize}
\item In \cite{Huke17_1} parts of the \textit{Ordering Phase} are
  tackled using two Mixed-Integer Linear Programs (MILPs).
\item  In  \cite{Hungerlaender2017492} we introduce the \textit{Slot Optimization Problem}.
 It describes the problem of determining the maximal number of available delivery time
windows for a new customer.
\item  Reference \cite{Hungerlaender201801} focuses on providing competitive MILP formulations
 for the \textit{Traveling Salesperson Problem with Time Windows}.

\end{itemize}
Next, we want to point out the most closely related literature:
\begin{itemize}
  \item
Reference \cite{camsa05}
describes a Home Delivery system that decides if a customer order is accepted.
Furthermore, the system assigns accepted orders to a time window under
consideration of the opportunity cost
of the orders.
In contrast to that, in our setup the customer takes the decision
to which delivery time window their order is assigned to.
\end{itemize}
The contribution of this work is given as followed:
\begin{itemize}
  \item A general description of an AHD system based on the use-case of an online
  grocery shopping service is given and a detailed description of the \textit{Ordering Phase} is provided.
  \item A heuristic solution approach for the introduced problem is presented. Further, the authors propose
  a hybrid approach that applies a Mixed-Integer Linear Program.
  \item Finally, novel benchmarks instances are introduced and computational experiments
  that evaluate the performance of the proposed approaches are presented.

\end{itemize}

\subsection{Challenges \& Key Steps of the Ordering Phase}

In this paper we focus on the \textit{Ordering Phase} and
 suggest solution approaches to
deal with the computational challenges arising during this phase.
In particular the runtime requirements for the optimization problems
applied during this phase are much more severe than in any other phase.

\begin{figure}[!ht]
\centering
  \begin{adjustbox}{width=0.8\textwidth}
    \begin{tikzpicture}[ every node/.style={font=\sffamily\footnotesize}]
     	\definecolor{blue}{rgb}{0.16, 0.32, 0.75}
    	\definecolor{red}{rgb}{0.7, 0.11, 0.11}
    	\definecolor{green}{rgb}{0.05, 0.5, 0.06}
    			\def \text{ 0.75cm}
    			\def \line {0.7mm}
          \def \linetwo {0.4mm}

    			\def \lvlA {0cm}
    			\def \lvlB {-0.5cm}
    			\def \lvlC {-1cm}
    			\def \lvlD {-1.5cm}
    			\def \color {darkgray}
    			\def \colorA{darkgray}
    			\def \scale {1}

    			\def \lvltext {1.2cm}
    			\draw [fill=lightgray!40, right] (-0.5,5) rectangle (10,0.3);

    			\node[fill=green!70,  line  width={\line}, draw=black, minimum width=1.9cm, minimum height=1cm,right]  at (0,1.2){\small 08:00-09:00};

    			\node[fill=red!70,  line  width={\line}, draw=black, minimum width=1.9cm, minimum height=1cm,right]  at (2.5,1.2){\small 09:00-10:00};
          \node[  line  width={\linetwo}, draw=black, minimum width=1.9cm, cross out, minimum height=1cm,right]  at (2.5,1.2){};

    			\node[fill=green!70,  line  width={\line}, draw=black, minimum width=1.9cm, minimum height=1cm,right]  at (5,1.2){\small 10:00-11:00};

    			\node[fill=green!70,  line  width={\line}, draw=black, minimum width=1.9cm, minimum height=1cm,right]  at (7.5,1.2){\small 11:00-12:00};

          \draw [fill=white, right] (-0.1,4.8) rectangle (9,4.4);
          \node[right] at(0,4.57)  {http://www.generic-grocery-store.com};

          \node[right] at(0,4)  {\large Order your groceries online!};

    	    \node[right] at(0,3.4)  {\textit{ Customer:} \quad John Doe};

    	    \node[right] at(0,3) {\textit{ Address:} \quad 123 Main St, Anytown, USA};

    	    \node[right] at(0,2.2) {\textit{Delivery Windows for:} \quad  March 3, 2018};

    			\end{tikzpicture}
\end{adjustbox}
  \caption{
  Illustration of a generic example website of an AHD service for grocery online shopping.
  Based on the customer's address, the system determines the
availability of the predefined  delivery time windows.
Non-available delivery time windows, e.g., 09:00-10:00, are crossed out.}
  \label{fig:exampleWebsite}
\end{figure}
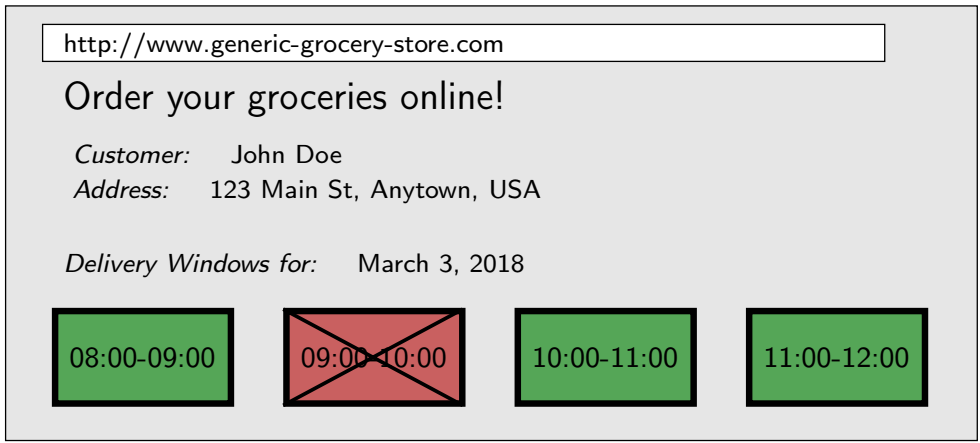

All decisions taken in the foregoing \textit{Tactical Planning Phase} are
 considered as input variables.
During the \textit{Ordering Phase} customers can book their grocery deliveries through a web service.
Figure  \ref{fig:exampleWebsite}
illustrates a generic example website and its main features.
Clearly, the web service should respond to the customer requests with as little
delay as possible.
Fetching and providing the input data for the booking process requires
communication across several services and many data base queries.
As this already requires a significant amount of time, there is even less
time to solve the actual optimization problem.

During the \textit{Ordering Phase}, the web service must
repeatedly solve an online variant of a \textit{capacitated Vehicle
  Routing Problem with Time Windows} (cVRPTW).
As the cVRPTW is known to be NP-hard \cite{lenstra1981}, the naive
approach of solving
a new cVRPTW instance from scratch for each new customer order is far from
being applicable in an online environment, even when using fast Meta-heuristics  \cite{el-sherbeny_2010},
due to the tight restrictions regarding runtime.

For clarity of exposition, we further split the \textit{Ordering
  Phase} into the following \textit{four key steps:}

\paragraph{Initialization step:}
The system sets up an empty delivery schedule, i.e., a cVRPTW instance
with a fixed number of vehicles and
corresponding operation times that were determined during the
\textit{Tactical Planning Phase},
but not having any customers assigned to yet.

\paragraph{Get TWs step - \textit{The system determines available delivery time windows}:}
Based on the current delivery schedule the system
determines which delivery time windows are available to a new 
customer. During times with high customer request rates this step has to
be performed within milliseconds.
The available time windows are then presented to the customer through the web service.
Note that the customer has to provide a delivery address such that
a routing system can estimate the travel times between all pairs of customers.

Optionally, for reasons of profit maximizing, some available time windows can be hidden
 from the customer or be offered at different rates.
 However, we do not consider any kind of slot pricing in this work.
For related and recent work on pricing in the context of AHD systems
we refer the reader to
\cite{yast14,Klein2017b,Klein2017}.

\paragraph{Set TW step - \textit{Customer books a delivery time window}:}
Using the website or app, a new customer selects her/his preferred delivery time window.
As other customers might have booked deliveries since the \textit{Get TWs} step,
the system must double-check if this delivery time window is still available.
If the insertion is still feasible,
 the system adds the new order into the
working schedule.
In order to avoid queuing issues during the critical \textit{Set TW}
step, the system does not allow any other
manipulations of the schedule.
\\
In case that the requested delivery time window
is not available anymore, the \textit{Get TWs} step is triggered again, and
an updated list of available delivery time windows is presented to the customer.

\paragraph{Improvement step:}
In this step the system tries to improve the working schedule
such that as many delivery time windows as possible can be offered
to potential future customers and therefore more customers can place
their orders.
Choosing the total travel time as
objective function has proven to be a reasonable
choice to achieve this goal.
While the fleet and the assignment of customers to time windows is fixed,
 the assignment of customers to  delivery vehicles as well as
 the sequences in which the vehicles visit the customers can be altered.
 Typically the \textit{Improvement} step may take several seconds, but during times with
 high customer request rates the step can be omitted or
triggered after, e.g.,\ every 10 \textit{Set TW} steps.

Note that at any time there is exactly one
 working schedule in the system.

In the following section we formally introduce the cVRPTW as the underlying
optimization problem of the \textit{Ordering
  Phase} and propose algorithmic strategies for dealing with the
cVRPTW during each of the four steps described above.

\section{Algorithms}
\label{sec:algorithms}
\subsection{Formal Definition of the cVRPTW}

In this section we now formally introduce the cVRPTW and some further required
notations.

\paragraph{Basic Definitions:}
A cVRPTW  instance is typically defined by the following input data:
\begin{itemize}
  \item A set of \textit{time windows} $\mathcal{W} =
    \left\{w_1,\dots,w_q
    \right\}$, where each time window  $w \in \mathcal{W}$ is
  defined through its \textit{start time} $s_w$ and its \textit{end
    time} $e_w$.
    We assume that the time windows are unique.
  	Hence, there do not exist time windows $w_a, w_b, \in
          \mathcal{W},  w_a \neq w_b$
  with $s_{w_a} = s_{w_b}$ and $e_{w_a} = e_{w_b}$.

\item A set of \textit{customers} $\mathcal{C}, ~|\mathcal{C}|=p$,
with corresponding \textit{order weight} function
$c: \mathcal{C} \to \mathbb{R}^{>0}$,
a \textit{service time} function $s: \mathcal{C} \to
\mathbb R^{>0}$,
and
a \textit{travel time} function $t: \mathcal{C} \times \mathcal{C} \to
\mathbb R^{>0}$ where we set the travel time from a customer $a$
to itself to $0$, i.e. $t(a,a)=0, ~a\in\mathcal{C}$.

\item A function $w: \mathcal{C} \to \mathcal{W}$ that
assigns to each customer a time window, during which the delivery vehicle
has to arrive at the customer.

\item A \textit{schedule}
$\mathcal{S}=\left\{\mathcal{A},\mathcal{B},\dots
\right\}$, consisting of $|\mathcal{S}|=m$ tours with assigned
  capacities $C_k,\ k \in \mathcal{S},$
where $C_k$ corresponds to the \textit{capacity} of the vehicle that
operates tour $k$.

\end{itemize}

A \textit{tour} $\mathcal{A} = \{ a_1, a_2,\dots, a_n\}$ contains $n$
customers, where the indices of the customers display the
sequence in which the customers are visited.
To improve clarity of exposition, we sometimes
additionally use upper indices, i.e.,  $\mathcal{A} = \{a_1^{w(a_1)},
a_2^{w(a_2)}, \dots, a_n^{w(a_n)}\}$, which indicate the time windows
assigned to the customers.
Furthermore, each tour $\mathcal{A}$ has assigned
\textit{start} and \textit{end times} that we denote as
{\it start$_{\mathcal{A}}$} and
{\it end$_{\mathcal{A}}$} respectively.
Hence, the vehicle executing tour $\mathcal{A}$ can leave from the start depot
no earlier than {\it start$_{\mathcal{A}}$}
and must return to the end depot
 no later than {\it end$_{\mathcal{A}}$}.

\paragraph{Structured Time Windows:}

Two  time windows $w_a$ and $w_b$ are \textit{non-overlapping} if and only if
$e_{w_a} \leq s_{w_b}$ or $ e_{w_b} \leq s_{w_a}$.
Therefore, $w_a$ and $w_b$  do overlap if and only if
$ s_{w_b} < e_{w_a}$ and $s_{w_a} < e_{w_b}$.
We speak of \textit{structured} time windows, if all time windows in $\mathcal{W}$
are pair-wise
non-overlapping and if the number of customers $|\mathcal{C}| = p$ is
much larger than the number of time windows
$|\mathcal{W}| = q$, i.e., $n \gg q$, and therefore typically
several customers
are assigned to the same time window.
We denote the corresponding variant of the cVRPTW as the
\textit{capacitated Vehicle Routing Problem with structured Time Windows} (cVRPsTW).
 Structured time windows are a specialty that arises
in the Attended Home Delivery use case, as well as some other modern routing
applications. Note that the corresponding assumptions do not impose severe restrictions
to the supplier nor the customers, but allow for a more efficient
optimization of the corresponding logistics operations.\\

\subsection{Arrival Times \& Feasibility}
Next let us give formal, recursive definitions of the earliest
and latest arrival times that are needed to define feasibility of a
schedule and of the insertion of a new customer.
\paragraph{Earliest and Latest Arrival Times:}
 We consider a fixed tour
 $\mathcal A =\{a_0, a_1, \dots,\allowbreak a_n, a_{n+1}\} $,
  where $a_0$ is the start
depot, $a_{n+1}$ is the end depot and $\{a_1,\ldots,a_{n}\}$ is the
set of customers assigned to tour $\mathcal A$. Note that all our approaches do not
move the depots. Hence, customers can only be inserted after the start depot
and before the end depot.
The earliest (latest) arrival time $\alpha_{a_i}$
($\beta_{a_i}$) gives the earliest (latest) time
at which the
vehicle may arrive at $a_i$, who is the
$i^{\text{th}}$ customer on the tour, while not violating
time window and travel time constraints on the preceding (subsequent)
tour:
\begin{align*}
		&\alpha_{a_0} :=start_{\mathcal A},
		&&\alpha_{a_{j+1}}:=
		\max
		\left\{
		s_{w(a_{j+1})},~
		\alpha_{a_j} + s(a_j)+t(a_j, a_{j+1})
		\right\},~
         j \in [n-1],
        \\
    &&&    \alpha_{a_{n+1}}:=
        \alpha_{a_n}+ s(a_n) +t(a_n, a_{n+1}),\\
		&\beta_{a_{n+1}}=end_{\mathcal A},
    &&\beta_{a_{j-1}}:=
		\min
		\left\{
		e_{w(a_{j-1})}, ~
		\beta_{a_j}-t(a_{j-1}  -s(a_{j-1}), a_{j})
		\right\}, \\
	&&&  j \in [n]\setminus \{1\}, \quad
        \beta_{a_0}:= \beta_{a_1} - t(a_0,a_1).
\end{align*}
\paragraph{Feasibility of a Schedule:}
Now we can concisely define the feasibility of a tour and a schedule with the
help of the earliest arrival times.
A schedule $\mathcal{S}$ is \textit{feasible}, if all its
tours are feasible. A tour $\mathcal{A}$ is \textit{feasible}, if
it satisfies both of the following conditions:
\begin{align*}
		s_{w(a_i)} \leq    \alpha_{a_i} &\leq e_{w(a_i)},   ~i \in [n],
    \quad \wedge \quad \alpha_{a_{n+1}}\leq end_{\mathcal{A}}
			& \qquad&
          \text{\meth{(TFEAS)}}, \\
		\sum\limits_{i \in [n]} c(a_{i})  &\leq C_{\mathcal{A}},  &\qquad
			&
        \text{\meth{(CFEAS)}}.
\end{align*}
While \meth{TFEAS} ensures that the arrival times at each customer are
within their assigned time windows,
\meth{CFEAS} ensures that the capacity of $\mathcal{A}$ is not exceeded.

\paragraph{Feasibility of an Insertion:}
Now we further use
the concepts of earliest and
 latest arrival time to facilitate and algorithmically speed up
 feasibility checks of tours
 after inserting an additional customer.
 A new customer $\tilde{a}^w$ can be feasibly inserted  \textit{with respect to
   time} between customers $a_i$ and $a_{i+1},\ i \in [n_0],$ into a
 feasible tour $\mathcal{A}$ if the following condition holds:
 \begin{align}
&&& \alpha_{\tilde{a}^w} \leq \beta_{\tilde{a}^w},\qquad
   \meth{TFEAS($\tilde{a}^w,{i+1},\mathcal{A}$)}, \label{cond:tfeas}
   \\
&\text{where,}
\quad
&&
\alpha_{\tilde{a}^w}:=
\max\{
s_w, ~\alpha_{a_i} + s(a_i) + t(a_i, \tilde{a}^w)
\},
 \nonumber\\
 &&&
 \beta_{\tilde{a}^w}:=
 \min \{
 e_w, ~\beta_{a_{i+1}} - s(\tilde{a}^w) - t(\tilde{a}^w, a_{i+1})
 \} .\nonumber
 \end{align}

\begin{figure}[!ht]
 \centering
  \begin{adjustbox}{width=0.9\textwidth}

    \tikzset{every picture/.style={line width=0.8pt}}

    \begin{tikzpicture}[every node/.style={font=\sffamily\footnotesize}]

    \draw[color=lightgray, line  width=1mm,-] (3cm,-1cm)  -- (3cm,1cm);
    \draw[color=lightgray,  line  width=1mm,-] (9cm ,-1cm)  -- (9cm,1cm);
    \node (text) at  (3cm, -1.25cm)  {$s_{w}$};
    \node (text) at  (9cm, -1.25cm)  {$e_{w}$};

    \node [	draw,xshift=1.47cm,yshift=0cm,	minimum width=1.1cm,
    minimum height=0.7cm,] (a_i) {$a_i^{w_{j-1}}$};
    \node [	draw,xshift=10cm,yshift=0cm,	minimum width=1.1cm,
    minimum height=0.7cm,] (a_i1) {$a_{i+1}^{w_{j+1}}$};
    \node [color=blue,	draw,xshift=7cm,yshift=0cm,	minimum
    width=0.8cm,		minimum height=0.7cm,] (a_w) {$\tilde{a}^{w_j}$};
    \node [color=blue,	draw,xshift=4.3cm,yshift=0cm,	minimum
    width=0.8cm,		minimum height=0.7cm,] (a_w_1)
    {$\tilde{a}^{w_j}$};

    \node (inv_2) at (11.0cm,0cm) {};
    \node (inv_1) at (0.5cm,0cm) {};

    \draw[draw, ->] (a_w.east) -> (a_i1.west)  ;
    \draw[draw, ->] (a_i.east) -> (a_w_1.west) ;
    \draw[decorate,decoration=snake] (inv_1.east) -> (a_i.west) ;
    \draw[decorate,decoration=snake] (a_i1.east) -> (inv_2.west) ;

    \draw[decorate,decoration=brace, color=black] ([shift={(0.02,0.05)}]
    a_i.east) -- ([shift={(-0.02,0.05)}] a_w_1.west)  node[above, midway ]
    {$\quad t(a_i^{w_{j-1}}, \tilde{a}^w)$ \hspace*{0.15cm} } (0,0);
    \draw[decorate,decoration=brace, color=black] ([shift={(0.02,0.05)}]
    a_w.east) -- ([shift={(-0.02,0.05)}] a_i1.west)  node[above, midway ]
    {$\quad t(\tilde{a}^{w_j}, a_{i+1}^{w_{j+1}})$ \hspace*{0.05cm} } (0,0);
    \draw[decorate,decoration=brace, color=black] ([shift={(0,0.05)}]
    a_w.north west) -- ([shift={(0,0.05)}] a_w.north east)  node[above,
    midway ] {$s(\tilde{a}^{w_{j}})$} (0,0);
    \draw[decorate,decoration=brace, color=black] ([shift={(0,0.05)}]
    a_i.north west) -- ([shift={(0,0.05)}] a_i.north east)  node[above,
    midway ] {$s(a_i^{w_{j-1}})$} (0,0);

    \draw[draw, ->] ([shift={(0,-0.7)}]a_w.south  west) --
    ([shift={(0,-0.05)}]a_w.south  west);
    \node (text) at ([shift={(0.4,-0.55)}] a_w_1.south west)  {$\alpha_{\tilde{a}^{w_j}}$};
    \draw[draw, ->] ([shift={(0,-0.7)}]a_w_1.south west) --
    ([shift={(0,-0.05)}]a_w_1.south west);
    \node (text) at ([shift={(0.4,-0.55)}]a_w.south  west)  {$\beta_{\tilde{a}^{w_k}}$};
    \draw[draw, ->] ([shift={(0,-0.7)}]a_i.south  west) --
    ([shift={(0,-0.05)}]a_i.south  west);
    \node (text) at ([shift={(0.55,-0.55)}] a_i.south west)  {$\alpha_{a_i^{w_{j-1}}}$};
    \draw[draw, ->] ([shift={(0,-0.7)}]a_i1.south  west) --
    ([shift={(0,-0.05)}]a_i1.south  west);
    \node (text) at ([shift={(0.55,-0.55)}] a_i1.south west)
    {$\beta_{a_{i+1}^{w_{j+1}}}$};

    \draw[decorate,decoration=brace, color=red]
    ([shift={(0,-0.75)}]a_w.south  west)
     --
    ([shift={(0,-0.75)}]a_w_1.south west)
      node[above, midway ]
    {$\leq$} (0,0);

    \draw[color= darkgray, draw,line  width=0.7mm, ->,>=stealth] (0.5cm,-1.6cm) -- (11cm, -1.6cm);
    \draw (6cm,-1.8cm) node {Time};

    \end{tikzpicture}

\end{adjustbox}
 \caption{Depiction of a feasible insertion with respect to time of $\tilde{a}^{w_j}$ between
   $a_i^{w_{j-1}}$ and $a_{i+1}^{w_{j+1}}$, i.e.,\
    \meth{TFEAS($\tilde{a}^{w_j},{i+1},\mathcal{A}$)} holds.   }
     \label{fig:feasInsert}
 \end{figure}
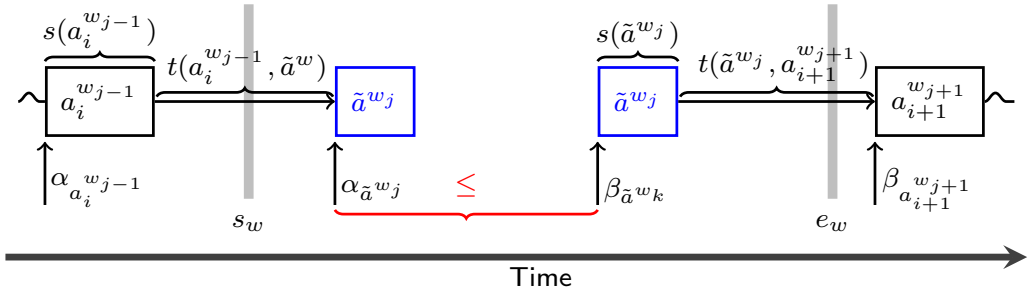

Condition \eqref{cond:tfeas} ensures that we arrive at customer $\tilde{a}^w$
 early enough, such that we can leave from $\tilde{a}^w$ early enough,
 to handle all subsequent customers of $\mathcal{A}$ within their
 assigned time windows. We refer to Figure \ref{fig:feasInsert} for an
 illustration of the above condition.

 Additionally, we have to check, if the sum of the weights of the
 customer orders assigned to tour $\mathcal{A}$, does not exceed the capacity
 $C_\mathcal{A}$. The insertion of $\tilde{a}^w$ into tour $\mathcal{A}$ is
   feasible  \textit{with respect to capacity}, if the following
   condition holds:
 \begin{align}
 \sum\limits_{i \in [n]} c(a_i)  + c(\tilde{a}^w) \leq C_\mathcal{A},\qquad
   \meth{CFEAS($\tilde{a},\mathcal{A}$)}. \label{cond:cfeas}
 \end{align}

Assuming that all earliest and latest arrival times and the sum of capacities
have already been calculated, Conditions \eqref{cond:tfeas} and \eqref{cond:cfeas}
allow to check the feasibility of an insertion of a new customer into
a given time window in $\bigO(1)$.

If we conduct an insertion that is feasible with respect to time and
capacity and decide to insert $\tilde{a}^w$ we receive a new tour
$\mathcal{\tilde{A}}=\{a_0,a_1,\dots,a_i,\tilde{a}^w,a_{i+1},\dots,\allowbreak a_n,a_{n+1}
\}$. Customer  $\tilde{a}^w$ is then assigned index $i+1$ and the indices of all
succeeding customers are increased by one.
Clearly, earliest and latest arrival times and the sum of capacities of the
modified tour must be updated, which can be done in $\bigO(n)$.

  Note that, in the context of an offline \textit{Traveling Salesperson Problem  with Time Windows},
   the Generalized Insertion Heuristic proposed in Reference
  \cite{gendreau_generalized_1998} uses concepts analog to our earliest
   and latest arrival times. The authors also check the feasibility of
   possible insertions with two
   conditions that resemble Condition \eqref{cond:tfeas}.
Due to their efficient computation,
Conditions \eqref{cond:tfeas} and \eqref{cond:cfeas} form the
basic building blocks of our Local Search heuristic
that we describe in the following subsections.
Moreover, the concept is flexible enough such that valuable extensions, such as
time-dependent travel times or the integration of driving breaks into the schedule,
can be employed without major changes.

\subsection{Local Search Heuristic}

We consider a Local Search heuristic that  uses two
neighborhoods for exchanging customer orders between two tours:
\begin{enumerate}
	\item The \textit{1-move} neighborhood moves a customer from one tour to
	another tour.
	\item  The \textit{1-swap} neighborhood swaps two customers between
	two different tours.
\end{enumerate}

Accordingly we define the \meth{1move}$(\tilde{a}^w,\mathcal{A},\mathcal{B})$
operation as the procedure where we \textit{remove} customer
$\tilde{a}^w$ from tour $\mathcal{A} \in \mathcal{S}$ and try to
feasible \textit{insert} it into
 tour $\mathcal{B} \in \mathcal{S}, ~\mathcal{A}\neq \mathcal{B},$
 within time window $w$.
If at least one feasible insertion position for $\tilde{a}$ in $\mathcal{B}$
is found that additionally decreases the total travel time of the delivery
schedule, we denote the \textit{1-move}
as \textit{improving}.

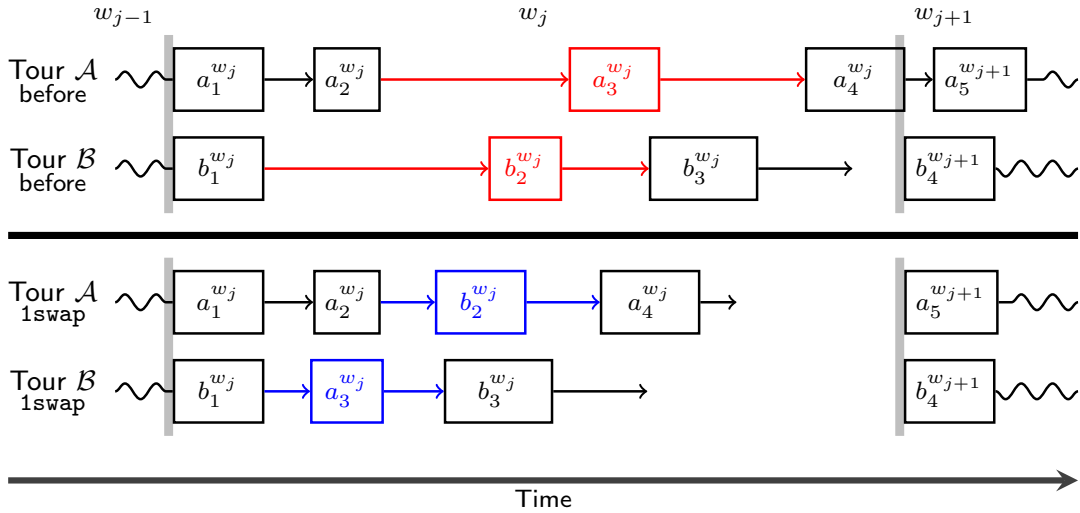
\begin{figure}[!ht]
	\centering
  \begin{adjustbox}{width=0.9\textwidth}
    \tikzset{every picture/.style={line width=0.8pt}}

    \begin{tikzpicture}[every node/.style={font=\sffamily\footnotesize}]

      \draw[color=lightgray, line  width=1mm,-] (10cm,-1cm)  -- (10cm,1cm);
      \draw[color=lightgray,  line  width=1mm,-] (1.8cm ,-1cm)  -- (1.8cm,1cm);
      \draw[color=black, line  width=0.8mm,-] (0cm, -1.25cm)  -- (12cm,-1.25cm);
      \draw[color=lightgray, line  width=1mm,-] (10, -1.5cm)  -- (10,-3.5cm);
      \draw[color=lightgray,  line  width=1mm,-] (1.8cm ,-1.5cm)  -- (1.8cm,-3.5cm);

      \draw  (1.3cm,1.2cm) node {$w_{j-1}$};
      \draw  (5.9cm,1.2cm) node {$w_j$};
      \draw (10.5cm,1.2cm) node {$w_{j+1}$};
    \draw [] (0.5cm,.5cm)  node {$ \stackrel{\text{\small{Tour $\mathcal{A}$}}} { \text{\footnotesize{before}}}$};
    \draw [] (0.5cm,-.5cm) node {$ \stackrel{\text{\small{Tour $\mathcal{B}$}}} { \text{\footnotesize{before}}}$};
    \draw [] (0.5cm,-2cm)  node {$ \stackrel{\text{\small{Tour $\mathcal{A}$}}} { \text{\footnotesize{\meth{1swap}}}}$};
    \draw [] (0.5cm,-3cm) node {$ \stackrel{\text{\small{Tour $\mathcal{B}$}}} { \text{\footnotesize{\meth{1swap}}}}$};

    \node [	draw,xshift=2.36cm,yshift=0.5cm,	minimum width=1cm,		minimum height=0.7cm,] (a2) {$a_1^{w_j}$};
    \node [	draw,xshift=3.8cm,yshift=0.5cm,	minimum width=0.6cm,		minimum height=0.7cm,] (a3) {$a_2^{w_j}$};
    \node [	draw,xshift=6.8cm,yshift=0.5cm,	minimum width=1cm,		minimum height=0.7cm, color=red] (a4) {$a_3^{w_j}$};
    \node [	draw,xshift=9.5cm,yshift=0.5cm,	minimum width=1.1cm,		minimum height=0.7cm,] (a5) {$a_4^{w_j}$};
    \node [	draw,xshift=10.9cm,yshift=0.5cm,	minimum width=0.8cm,		minimum height=0.7cm,] (a6) {$a_5^{w_{j+1}}$};
    \draw[draw, ->] (a2.east) -> (a3.west) ;
    \draw[draw,color=red, ->] (a3.east) -> (a4.west) ;
    \draw[draw,color=red, ->] (a4.east) -> (a5.west) ;
    \draw[draw, ->] (a5.east) -> (a6.west) ;
    \draw[decorate,decoration=snake] (1.2cm,0.5cm) -> (a2.west) ;
    \draw[decorate,decoration=snake] (12cm,0.5cm) -> (a6.east) ;

    \node [	draw,xshift=2.36cm,yshift=-0.5cm,	minimum width=1cm,		minimum height=0.7cm,] (b2) {$b_1^{w_j}$};
    \node [draw,xshift=5.8cm,yshift=-0.5cm,	minimum width=0.8cm,	minimum height=0.7cm,, color=red] (b3) {$b_2^{w_j}$};
    \node [	draw,xshift=7.8cm,yshift=-0.5cm,	minimum width=1.2cm,		minimum height=0.7cm] (b4) {$b_3^{w_j}$};
    \node (invB) at (9.6cm,-0.5cm) {};
    \node [draw,xshift=10.56cm,yshift=-0.5cm,	minimum width=1cm,	minimum height=0.7cm,] (b5) {$b_4^{w_{j+1}}$};
    \draw[draw, color=red,->] (b2.east) -> (b3.west) ;
    \draw[draw,color=red, ->] (b3.east) -> (b4.west) ;
    \draw[draw, ->] (b4.east) -> (invB.west) ;
    \draw[decorate,decoration=snake] (1.2cm,-0.5cm) -> (b2.west) ;
    \draw[decorate,decoration=snake] (12cm,-0.5cm) -> (b5.east) ;

    \node [	draw,xshift=2.36cm,yshift=-2cm,	minimum width=1cm,		minimum height=0.7cm,] (a2) {$a_1^{w_j}$};
    \node [	draw,xshift=3.8cm,yshift=-2cm,	minimum width=0.6cm,		minimum height=0.7cm,] (a3) {$a_2^{w_j}$};
    \node [	draw,xshift=5.3cm,yshift=-2cm,	minimum width=1cm,		minimum height=0.7cm, color=blue] (a4) {$b_2^{w_j}$};
    \node [	draw,xshift=7.2cm,yshift=-2cm,	minimum width=1.1cm,		minimum height=0.7cm,] (a5) {$a_4^{w_j}$};
    \node [	draw,xshift=10.58cm,yshift=-2cm,	minimum width=0.8cm,		minimum height=0.7cm,] (a6) {$a_5^{w_{j+1}}$};
    \draw[draw, ->] (a2.east) -> (a3.west) ;
    \draw[draw, ->, color=blue] (a3.east) -> (a4.west) ;
    \draw[draw, ->, color=blue] (a4.east) -> (a5.west) ;
    \node (inv) at (8.3cm,-2cm ) {};

    \draw[draw, ->] (a5.east) -> (inv.west) ;
    \draw[decorate,decoration=snake] (1.2cm,-2cm) -> (a2.west) ;
    \draw[decorate,decoration=snake] (12cm,-2cm) -> (a6.east) ;

    \node [	draw,xshift=2.36cm,yshift=-3cm,	minimum width=1cm,		minimum height=0.7cm,] (b2) {$b_1^{w_j}$};
    \node [draw,xshift=3.8cm,yshift=-3cm,	minimum width=0.8cm,	minimum height=0.7cm,color=blue] (b3) {$a_3^{w_j}$};
    \node [	draw,xshift=5.5cm,yshift=-3cm,	minimum width=1.2cm,		minimum height=0.7cm] (b4) {$b_3^{w_j}$};
    \node (invB) at (7.3cm,-3cm) {};
    \node [draw,xshift=10.56cm,yshift=-3cm,	minimum width=1cm,	minimum height=0.7cm,] (b5) {$b_4^{w_{j+1}}$};
    \draw[draw, ->, color=blue] (b2.east) -> (b3.west) ;
    \draw[draw, ->, color=blue] (b3.east) -> (b4.west) ;
    \draw[draw, ->] (b4.east) -> (invB.west) ;
    \draw[decorate,decoration=snake] (1.2cm,-3cm) -> (b2.west) ;
    \draw[decorate,decoration=snake] (12cm,-3cm) -> (b5.east) ;

    \def \lvlTimeLine {-4cm}
    \draw[color= darkgray, draw,line  width=0.7mm, ->,>=stealth] (0cm,{\lvlTimeLine}) -- (12cm, {\lvlTimeLine});
    \draw (6cm,{\lvlTimeLine -0.2cm}) node {Time};
    \end{tikzpicture}\end{adjustbox}
\caption{
Reduction of the total travel time of a schedule induced by an improving
\meth{1swap}$(a_{3}^{w_j}, \mathcal{A}, \mathcal{B})$ operation.
}
	\label{fig:oneswap}
\end{figure}

As a  \meth{1swap}$(\tilde{a}^w,\mathcal{A},\mathcal{B})$
operation we define the procedure where we
try to \textit{exchange} customer $\tilde{a}^w$
with any customer within assigned time window $w$ from a different tour
$\mathcal{B}$.
If at least one such exchange
  decreases  the total travel time
of the schedule, we denote the
\textit{1-swap} as \textit{improving}.
In general we always select the exchange of an improving \textit{1-swap}
that results in the largest decrease of the total travel time of the
delivery schedule.
In Figure \ref{fig:oneswap} we provide an illustration of an improving
\meth{1swap}$(a_{3}^{w_j}, \mathcal{A}, \mathcal{B})$.

\subsection{Algorithmic Strategies}
\label{subsec:algStrategies}

In this subsection we describe how to combine
the Local Search heuristics presented
in the previous subsection in order to conduct sufficiently fast
\textit{Get TWs}, \textit{Set TW} and  \textit{Improvement} steps.

\paragraph{The \textit{Get TWs} step:}
In this step we aim to quickly identify all time windows during
which a new customer $\tilde{a}$ can be inserted into (at least one
of) the current tours.
We suggest to use the following procedure:
\begin{itemize}
\item \textit{Simple Insertion:}
For each time window $w \in \mathcal{W}$ iterate over all tours
$\mathcal{A} \in \mathcal{S}$ and all possible insertion points within $w$
and check Conditions \eqref{cond:tfeas} and  \eqref{cond:cfeas}.
A  time window $w$ is considered as being available, if both conditions hold
for at least one insertion point. In this case we add $w$ to the set
of available time windows $\mathcal{T}_{\tilde{a}} \subseteq \mathcal{W}$.
This procedure is  computationally very cheap
	and runs
  within 1 millisecond (ms) for all benchmark instances considered
  in our computational study.
\end{itemize}
Then the time windows $\mathcal{T}_{\tilde{a}}$ are offered to the
customers through the web service.
Note that Reference \cite{camsa05} proposes a similar procedure for the VRP,
i.e., without considering time windows.

Additionally to this \textit{Simple Insertion} heuristic we
introduced an Adaptive Neighborhood Search (ANS) in
\cite{Hungerlaender2017492} that is especially tailored to the
\textit{Get TWs} step.
The ANS applies \textit{1-move} and \textit{1-swap}
operations to free up time during a specific time window on a selected
tour in order to enable the insertion of the new customer.
ANS has proven to find more available time windows than \textit{Simple
  Insertion}, while still being fast enough for most applications, as
long as the customer request rate is moderate.

\paragraph{The \textit{Set TW} step:}
Once the customer has selected a time window $\tilde{w}$ from the set
$\mathcal{T}_{\tilde{a}} $,
we double check its availability in the same manner as in the
\textit{Get TWs} step, and then we immediately insert $\tilde{a}$
 into $\tilde{w}$
at the best found insertion point.

\paragraph{The Improvement step:}

In this step we aim to reduce the
total travel time
 of the delivery schedule
by using one of the
following two procedures:
\begin{itemize}
	\item  \textbf{Local-improvement}: Our computationally
          cheap, yet quite effective Local Search heuristic
          builds the foundation of the \textit{Improvement} step.
          We combine \textit{1-move} and
				\textit{1-swap} operations, where we
                                focus on the \textit{1-move}
                                operations
				when possible, as
                                they are computationally cheaper and
                                in general more
                                effective
				than the \textit{1-swap} operations. We
                                stop our Local Search heuristic
				once we reach a local minimum of our
                                objective function with respect to
				our neighborhoods.
	\item \textbf{Local+TSPTW-improvement}: After the Local Search
          heuristic we additionally use MILPs
          proposed in our previous paper \cite{Hungerlaender201801}
          for optimizing all single tours
           that have changed since the last improvement step. In  \cite{Hungerlaender201801} we
                                         motivated and analyzed
the \textit{Traveling Salesperson Problem with Time
  Windows} (TSPTW) that
is a subproblem of the cVRPTW as each tour of the delivery schedule
corresponds to a TSPTW instance. Optimizing the single  tours of a
schedule to optimality has  been proven
to be critical to ensure driver satisfaction.
Hence, it ensures that drivers  not encounter any obvious loops on
their routes \cite{doi:10.1287/inte.2016.0875}. Also note that we use the current tours of
                                         our delivery schedule for
                                         warm starting the TSPTW
                                         MILPs.
	\end{itemize}
During the \textit{Ordering Phase} our Local Search heuristic
only performs improving operations. However, the algorithms can be simply
altered into a \textit{Simulated Annealing} approach by allowing also non-improving
operations, which is more appropriate for the \textit{Preparation
  Phase} when there is more time available for optimization.

\section{Computational Experiments}
\label{sec:experiments}
In this section, we present computational results on a set of benchmark
instances that are motivated by an online  shopping service of
an international grocery retailer.
We restrict our experiments to
a setup with structured time windows
as it has proven to be much more computationally efficient than
having overlapping time windows
and therefore is better suited for the use in an AHD system.
Accordingly
we consider the \textit{Traveling Salesperson Problem with structured Time
  Windows} (TSPsTW)
\cite{Hungerlaender201801},
 which is a special case of the cVRPsTW as each
tour of
the delivery schedule corresponds to a TSPsTW instance.

\subsection{Benchmark Instances}

In order to provide meaningful computational experiments we created
a benchmark set that resembles
real-world data focusing on urban settlement structures.
Moreover, the  benchmark instances are designed to reflect instances
as they arise in an online grocery shopping service
of a major international supermarket chain,
regarding
travel times, length of time windows, duration of service
times, customer order weights and their proportions to vehicle
capacities.
All instances can be downloaded from
\url{http://tinyurl.com/vrpstw}.
Note that the well-known VRPTW benchmark instances proposed by Solomon
\cite{solomon_algorithms_1987}
do not comply with our considered use-case.

In more detail our benchmark instances have the following characteristics:
\begin{itemize}

\item
\textbf{Grid Size.}
We consider a $20$\,km $\times$ $20$\,km square grid,
which is  roughly  of the size of Vienna as well as
a smaller grid of size $10$ \,km $\times$ $10$\,km that corresponds to smaller cities.

\item \textbf{Placing of Customers.}
In order to achieve varying customer densities,
only 20\,\% of the customer locations have been sampled from a two-dimensional
uniform distribution.
The remaining 80\,\% of the customer locations have been
randomly assigned to $10$ clusters.
The centers resp. shapes of those clusters have been randomly sampled as well.
Finally, the customer locations have been sampled from the assigned cluster.

\item \textbf{Depot location.}
We consider two different placements of the depot:
 at the center of the grid, and at  the center of the top left quadrant.
In each test setup there are equally many instances for both variants.
\item \textbf{Travel speeds.}
As proposed by Reference \cite{pan2017}
 we assume a travel speed of $20$\,km/h.
This number can be further supported by a recent report by
Vienna Public Transport
 \cite{WienerLinienEN2018}, where an average
travel speed for their fleet of buses of $17.7$\,km/h  during the day, $17.2$\,km/h at peak times,
and $20.0$\,km/h during evening
hours has been reported.

For the sake of simplicity, the distance between two locations is
 calculated as the Euclidean distance between them.
Travel times are calculated proportional to the Euclidean distances, using the assumed
travel speed.

 \item \textbf{Order  weights.} The order weights of customers have been
  sampled from a truncated normal distribution with
mean of 7 and standard deviation of 2, where the lower bound is 1 and
the upper bound is 15.

\item \textbf{Customer choice model.}
A customer choice model simulates the decisions that are usually taken by the
customers.
We choose a  simple model, where every customer has just one desired delivery time window
that has been set beforehand in the benchmark instance.
If the preferred time window is not offered to the customer, we assume that
the customer does not place an order and leaves the website.
We simulate this by a random assignment, following a uniform distribution,
of each customer to one time window out of a set of $5$ resp. $10$ consecutive time windows,
 where each is one hour long.
The later reflects a situation where customers choose from one-hour delivery time windows
 between 08:00 and 18:00.
Note that this is in contrast to real-world applications, where usually
certain time windows are more prominent among customers than others.
However, we chose a uniform distribution to obtain unbiased
results that allow for an easier
identification and clearer interpretation of the key findings.
\item \textbf{Service times.}
We assume the service time at each customer to be 5 minutes.
\end{itemize}

Typically, in the online grocery delivery use-case no more than 500 customer
are served from the same depot on a given day. Hence,
for our computational experiments we consider two benchmark sets that contain 500 customers each:
  \begin{enumerate}
  \item \label{setup:many}
  A benchmark set with {\it many short}  (30-40) tours
  having a capacity of 100 and 5 time windows each.
  \item \label{setup:moderate}
  A benchmark set with {\it fewer  long} (10-20) tours
  having a capacity of 200 and 10 time windows each.
  \end{enumerate}
Moreover, we simulate sparsely such as densely populated delivery regions
by using grids of size  $20$\,km $\times$ $20$\,km  resp.  $10$\,km $\times$ $10$\,km.

  \subsection{Experimental Setup}\label{ssec:setup}

  All experiments were performed on an Ubuntu 14.04 machine equipped
   with an Intel Xeon E5-2630V3 @ 2.4 GHz 8 core processor and 132 GB RAM.
  We implemented all algorithms in Java
  version 8 and use
  Gurobi 8.0.1 as IP-solver in single thread mode.
  We compare
  the algorithmic strategies presented in the previous section for both
  the \textit{Get TWs} and the \textit{Improvement} step.

  In all our experiments we iteratively insert new customers into the
  schedule, simulating
  customers placing orders online following the customer choice model.
  Hence, we assume that
  if the preferred time
  window is not offered to the customer he or she refuses to place an order and hence, the
  customer is not inserted into the schedule.
  Due to the iterative setup
  we can omit
   the \textit{Set TW} step and insert the new customer without
  double-checking the availability of the selected delivery time slot.
  We determine the following metrics averaged over 100 instances each:
  \begin{itemize}
  \item \textit{Get TWs} step:
  \begin{itemize}
  \item  Average  number of feasible time windows determined for each customer:
         corresponds to the number of time windows in which the order
         can be inserted.

  \item Average runtime of the \textit{Get TWs} step.
  \end{itemize}

  \item \textit{Improvement} step:
  \begin{itemize}
    \item Average reduction of the
      sum of travel times over all tours (given as percentage).

    \item Average reduction of the sum of travel times over all tours relative to the
    increase of travel time caused by the insertion of the new customer
    (given as percentage).
    \item Average number of TSPsTW MILPs solved.
    \item Average runtime of each \textit{Improvement} step.
  \end{itemize}
  \end{itemize}

\subsection{Results}\label{ssec:results}

Now let us present the results of our computational evaluation.
We examine the performance of our approaches on  instances with 500 customers.
The results for the sparse resp.  dense  benchmark sets
 are summarized in Table~\ref{tab:insertion_sparse}  resp.
Table~\ref{tab:insertion_dense} for the  \textit{Get TWs} step, and in
Table~\ref{tab:improvement_sparse} resp.
Table~\ref{tab:improvement_dense}  for the \textit{Improvement} step.

\begin{table}[ht!]
  \caption{Results for the \textit{Get TWs} step for our sparse benchmark scenarios.}
  \begin{center}
    \begin{tabular}{|l|rrr|rrr|}
    \hline
    {\bf  Get TWs step} & \multicolumn{3}{c|}{500 customers} &
    \multicolumn{3}{c|}{500 customers} \\
    Simple-insertion
    & \multicolumn{3}{c|}{100 capacity units } & \multicolumn{3}{c|}{200 capacity units} \\
    & \multicolumn{3}{c|}{5 time windows}
      & \multicolumn{3}{c|}{10 time windows}  \\
    Tours: & 30 & 35 & 40 & 10 & 15 & 24 \\
    \hline
    \hline
    {\bf  Average runtime} (sec:ms)  &
    0:001  & 0:001  & 0:001  &
    0:001  & 0:001  & 0:001    \cr
    \hline
    {\bf Number of  time }  &&&&&&\\
    {\bf windows offered} (avg.)
     & 4.29 & 4.93 & 5.00 &
    5.76 & 8.59 & 10.00 \cr \hline
      {\bf Total customers}  &&&&&&\\
    {\bf inserted } (avg.)
    & 428.9 & 493.4 & 500 &
    288.3 & 429.2 & 500 \cr
    \hline
    \end{tabular}
\end{center}
  \label{tab:insertion_sparse}
\end{table}

\begin{table}[ht!]
  \caption{Results for the \textit{Get TWs} step for our dense benchmark scenarios.}
\begin{center}
  \begin{tabular}{|l|rrr|rrr|}
  \hline
  {\bf  Get TWs step} & \multicolumn{3}{c|}{500 customers} &
  \multicolumn{3}{c|}{500 customers} \\
  Simple-insertion
  & \multicolumn{3}{c|}{100 capacity units } & \multicolumn{3}{c|}{200 capacity units} \\
  & \multicolumn{3}{c|}{5 time windows}
    & \multicolumn{3}{c|}{10 time windows}  \\
  Tours: & 30 & 35 & 40 & 10 & 15 & 20 \\
  \hline
  \hline
  {\bf  Average runtime} (sec:ms)  &
  0:001  & 0:001  & 0:001  &
  0:001  & 0:001  & 0:001  \cr

  \hline
  {\bf Number of  time }  &&&&&&\\
  {\bf windows offered} (avg.)
   & 4.28 & 4.94 & 5.00 &
  5.76 & 8.59 & 10.00 \cr \hline
  {\bf Total customers}  &&&&&&\\
  {\bf inserted} (avg.)
  & 428.3 & 493.7 & 500
   & 287.8 & 429.5 & 500\cr
  \hline
  \end{tabular}
\end{center}
  \label{tab:insertion_dense}
\end{table}

\begin{table}[ht!]
  \caption{Results for the \textit{Improvement} step for our sparse benchmark scenarios.}
  \begin{center}

    \begin{tabular}{|l|rrr|rrr|}
    \hline
    {\bf  Improvement step }  & \multicolumn{3}{c|}{500 customers} &
    \multicolumn{3}{c|}{500 customers} \\
    & \multicolumn{3}{c|}{100 capacity units } & \multicolumn{3}{c|}{200 capacity units} \\
    & \multicolumn{3}{c|}{5 time windows}
      & \multicolumn{3}{c|}{10 time windows}  \\
      Tours: & 30 & 35 & 40 & 10 & 15 & 20 \\
    \hline
    \hline
    {\bf  Avg. runtime (sec:ms)} &&&&&&\\
    Local &
    1:024 & 1:453 & 1:504 &
    0:653 & 1:512 & 1:966\cr
    Local+TSPsTW &
    1:115   & 1:547 & 1:600 &
    0:729 & 1:602 & 2:058\cr
    \hline
    {\bf Avg. improvement} &&&&&&\\
    {\bf over insertion step (\%)} &&&&&&\\
    Local &
    0.96 & 0.90 & 0.90&
    0.71 & 0.65 & 0.60\cr
    Local+TSPsTW &
    1.00  & 0.94 & 0.93&
    0.78 & 0.72 & 0.65\cr
    \hline
    {\bf Avg. improvement of} &&&&&&\\
    {\bf cost of insertion (\%)} &&&&&&\\
    Local &
    66.12 & 67.97 & 68.49&
    50.13 & 56.73 & 58.96\cr
    Local+TSPsTW &
    68.81 & 78.55 & 71.05&
    55.16 & 61.27 & 63.42\cr
    \hline
    {\bf  Avg. number of} &&&&&&\\
    {\bf  TSPsTW MILPs solved} &&&&&&\\
    Local+TSPsTW &
    3.48 & 3.76 & 3.81&
    1.82 & 2.25 & 2.41\cr
    \hline
    \end{tabular}
\end{center}
  \label{tab:improvement_sparse}
\end{table}

\begin{table}[ht!]
  \caption{Results for the \textit{Improvement} step for our dense benchmark scenarios.}
  \begin{center}
    \begin{tabular}{|l|rrr|rrr|}
    \hline
    {\bf  Improvement step }  & \multicolumn{3}{c|}{500 customers} &
    \multicolumn{3}{c|}{500 customers} \\
    & \multicolumn{3}{c|}{100 capacity units } & \multicolumn{3}{c|}{200 capacity units} \\
    & \multicolumn{3}{c|}{5 time windows}
      & \multicolumn{3}{c|}{10 time windows}  \\
    Tours: & 30 & 35 & 40 & 10 & 15 & 20 \\
    \hline
    \hline
    {\bf  Avg. runtime (sec:ms)} &&&&&&\\
    Local &
    1:468 & 1:882 & 1:912 &
    0:807 & 2:035 & 2:707\cr
    Local+TSPsTW &
    1:564 & 1:992 & 2:014&
    0:886 & 2:140 & 2:807\cr
    \hline
    {\bf Avg. improvement} &&&&&&\\
    {\bf over insertion step (\%)} &&&&&&\\
    Local &
    1.00 & 0.95 & 0.95&
    0.68 & 0.65 & 0.62\cr
    Local+TSPsTW &
    1.05 & 0.99 & 0.98 &
    0.76 & 0.70 & 0.66\cr
    \hline
    {\bf Avg. improvement of } &&&&&&\\
    {\bf  cost of insertion (\%)} &&&&&&\\
    Local &
    67.37 & 69.59 & 70.20&
    49.01 & 57.18 & 60.65\cr
    Local+TSPsTW &
    70.43 & 72.48 & 73.09&
    54.56 & 62.23 & 65.66\cr
    \hline
    {\bf  Avg. number of} &&&&&&\\
    {\bf  TSPsTW MILPs solved} &&&&&&\\
    Local+TSPsTW &
    3.59 & 3.89 & 3.94&
    1.93 & 2.37 & 2.56\cr
    \hline
    \end{tabular}
\end{center}
  \label{tab:improvement_dense}
\end{table}

The first benchmark set with many short tours corresponds to the left
column, and the results for the second benchmark set with few long tours
are displayed in the right column.
First, we observe that
the runtimes for both the  \textit{Get TWs} and the \textit{Improvement} step
are very low despite of the
large instances, which demonstrates that our solving approaches scale
very well.
It is worth pointing out that the  \textit{Get TWs} step stays below 1 ms
even for a large number of customers. This is crucial in order to deal with
high customer request rates at peak times.
Considering that between two \textit{Improvement} steps the schedule is altered
only by insertion of one customer, a reduction of our objective
function by $0.60\,\%$ to $1.05\,\%$
per step is remarkable.
This can be further underlined by the reported average
reduction of the cost of inserting the new customer which ranges from
$50.13\,\%$ to $78.55\,\%$.
Furthermore, we notice a moderate improvement of the hybrid heuristics over the
Local-improvement heuristics.

When comparing dense to sparse instances, we notice that the average
runtimes for the  \textit{Improvement} step are significantly higher on dense instances.
Further, we notice higher average runtimes on the benchmark set with fewer long tours.
The reason for both observations lies in the
larger number of customers per tour which makes these
instances more difficult to solve.
We also
observe that the reduction of the objective function achieved by the hybrid heuristics,
compared to the Local-improvement heuristics, is similar on both benchmark sets.

In summary, our suggested algorithms perform very good an both
benchmark sets as they are able to produce
delivery
schedules on large-scale instances
within the tight runtime restrictions imposed by the considered
application.

\section{Conclusion}
\label{sec:conclusion}
In this work, we considered an Attended Home Delivery (AHD) system in the context of  an
 online grocery shopping service offered by an international grocery retailer.
AHD systems are used whenever
a supplying company offers online shopping services that require that customers must be
present when their deliveries arrive.
 Therefore, the supplying company and the
customer must both agree on a  time window, which ideally is rather short,
 during which delivery is guaranteed.

In this paper we considered the overall fulfillment process of the AHD
system that can be described by four consecutive phases:
(1) Tactical Planning, (2) Ordering, (3) Preparation, and (4) Delivery.
We focused on
 the ordering phase, during which customers place their orders through the
 web service.
 Generally, this phase is the most challenging phase of an AHD system
 from a computational point of view.
 As for most AHD approaches in the literature, we considered
 a capacitated Vehicle Routing Problem with
 Time Windows  as the underlying optimization problem of the ordering phase.
 The online characteristic of this phase requires that the delivery schedule is
  built dynamically as  new orders are placed.
We split the computations of the ordering phase into four key steps
and proposed a  solution approach that
 allows to
 (non-stochastically) determine which delivery time windows can be
offered to potential customers.
Furthermore we
employed
a Local Search heuristic
to improve the delivery schedule and we also suggested a hybrid
approach that additionally to the Local Search heuristic employs
MILPs, which optimize single tours.

Finally, in an experimental evaluation, we demonstrated the efficiency of our
approaches on benchmark sets that are motivated by
an online grocery shopping service.
We considered  the capacitated Vehicle Routing Problem with
 structured Time Windows (cVRPsTW) for our benchmarking experiments.
 The special feature of the
 cVRPsTW is the additional structure of the time windows
which does not impose severe restrictions neither to the supplying company nor
to the customers.
Our computational study showed that the suggested algorithms can solve the
considered
cVRPsTW instances fast
enough to comply with the very strict runtime restrictions as they arise in AHD systems with
high customer request rates.



\begin{thebibliography}{10}
\providecommand{\url}[1]{\texttt{#1}}
\providecommand{\urlprefix}{URL }
\providecommand{\doi}[1]{https://doi.org/#1}

\bibitem{agca08}
Agatz, N., Campbell, A.M., Fleischmann, M., Savelsbergh, M.W.P.: Challenges and
  opportunities in attended home delivery. In: Golden, B., Raghavan, S., Wasil,
  E. (eds.) The Vehicle Routing Problem: Latest Advances and New Challenges,
  pp. 379--396. Springer US (2008)

\bibitem{camsa05}
Campbell, A.M., Savelsbergh, M.W.P.: Decision {S}upport for {C}onsumer {D}irect
  {G}rocery {I}nitiatives. Transportation Science  \textbf{39}(3),  313--327
  (2005)

\bibitem{Ehmke2012}
Ehmke, J.F.: Attended home delivery. In: Integration of Information and
  Optimization Models for Routing in City Logistics, pp. 23--33. Springer US,
  Boston, MA (2012). \doi{10.1007/978-1-4614-3628-7\_3}

\bibitem{Ehmke2014}
Ehmke, J.F., Campbell, A.M.: {Customer acceptance mechanisms for home
  deliveries in metropolitan areas}. European Journal of Operational Research
  \textbf{233}(1),  193--207 (2014).
  \doi{https://doi.org/10.1016/j.ejor.2013.08.028}

\bibitem{el-sherbeny_2010}
El-Sherbeny, N.A.: {Vehicle routing with time windows: {An} overview of exact
  heuristic and metaheuristic methods}. Journal of King Saud University
  \textbf{22},  123--131 (2010)

\bibitem{gendreau_generalized_1998}
Gendreau, M., Hertz, A., Laporte, G., Stan, M.: A {Generalized} {Insertion}
  {Heuristic} for the {Traveling} {Salesman} {Problem} with {Time} {Windows}.
  Operations Research  \textbf{46}(3),  330--335 (1998)

\bibitem{doi:10.1287/inte.2016.0875}
Holland, C., Levis, J., Nuggehalli, R., Santilli, B., Winters, J.: {UPS
  Optimizes Delivery Routes}. Interfaces  \textbf{47}(1),  8--23 (2017).
  \doi{10.1287/inte.2016.0875}

\bibitem{Huke17_1}
Hungerl\"ander, P., Maier, K., P\"ocher, J., Rendl, A., Truden, C.: Solving an
  on-line capacitated vehicle routing problem with structured time windows. In:
  Fink, A., F\"ugenschuh, A., Geiger, M.J. (eds.) Operations Research
  Proceedings 2016, pp. 127--132. Springer International Publishing, Basel
  (2017). \doi{10.1007/978-3-319-55702-1\_18}

\bibitem{Hungerlaender2017492}
Hungerl\"ander, P., Rendl, A., Truden, C.: {On the Slot Optimization Problem in
  On-Line Vehicle Routing }. Transportation Research Procedia  \textbf{27},
  492--499 (2017). \doi{10.1016/j.trpro.2017.12.046}

\bibitem{Hungerlaender201801}
Hungerl\"ander, P., Truden, C.: {Efficient and Easy-to-Implement Mixed-Integer
  Linear Programs for the Traveling Salesperson Problem with Time Windows}.
  Transportation Research Procedia  (2018),
  \url{http://www.optimization-online.org/DB_HTML/2018/02/6466.html}, accepted

\bibitem{Klein2017b}
Klein, R., Mackert, J., Neugebauer, M., Steinhardt, C.: {A model-based
  approximation of opportunity cost for dynamic pricing in attended home
  delivery}. OR Spectrum  (Dec 2017). \doi{10.1007/s00291-017-0501-3}

\bibitem{Klein2017}
Klein, R., Neugebauer, M., Ratkovitch, D., Steinhardt, C.: {Differentiated Time
  Slot Pricing Under Routing Considerations in Attended Home Delivery}.
  Transportation Science  (2017). \doi{10.1287/trsc.2017.0738}

\bibitem{lenstra1981}
Lenstra, J.K., Kan, A.H.G.R.: {Complexity of vehicle routing and scheduling
  problems}. Networks  \textbf{11}(2),  221--227 (1981).
  \doi{10.1002/net.3230110211}

\bibitem{pan2017}
Pan, S., Giannikas, V., Han, Y., Grover-Silva, E., Qiao, B.: {Using
  customer-related data to enhance e-grocery home delivery}. Industrial
  Management \& Data Systems  \textbf{117}(9),  1917--1933 (2017).
  \doi{10.1108/IMDS-10-2016-0432}

\bibitem{solomon_algorithms_1987}
Solomon, M.M.: Algorithms for the {Vehicle} {Routing} and {Scheduling}
  {Problems} with {Time} {Window} {Constraints}. Operations Research
  \textbf{35}(2),  254--265 (1987). \doi{10.1287/opre.35.2.254},
  \url{https://doi.org/10.1287/opre.35.2.254}

\bibitem{Syndy2015}
Syndy: {The state of online grocery retail in Europe 2015} (2015),
  \url{http://www.syndy.com/report-the-state-of-online-grocery-retail-2015/}

\bibitem{WienerLinienEN2018}
{Vienna Public Transport (Wiener Linien)}: {2017 - Facts and Figures} (2018),
  \url{https://www.wienerlinien.at/media/files/2018/facts_and_figures_2017_243486.pdf}

\bibitem{yast14}
Yang, X., Strauss, A.K., Currie, C.S.M., Eglese, R.: Choice-{B}ased {D}emand
  {M}anagement and {V}ehicle {R}outing in {E}-{F}ulfillment. Transportation
  Science  \textbf{50}(2),  473--488 (2016)

\end{thebibliography}
\end{document}